\documentclass[a4paper]{article}  
\usepackage[latin1]{inputenc}
\usepackage[T1]{fontenc}
\usepackage[english]{babel}
\usepackage{amssymb,amsmath,mathtools,bm}

\usepackage[margin=3cm]{geometry}
\usepackage{calc}

\usepackage[all,line,dvips,cmtip]{xy}
\CompileMatrices

\newcommand\cxymatrix[2][]{%
  \begin{xy}*[c]\xybox{\xymatrix#1{#2}}\end{xy} }

\makeatletter

\newcommand*\DeclareMathSymbolShorthand[2]{
   \begingroup
   \setkeys{DMSS}{name=#2,#1}%
   \if\DMSS@overwrite 
   \else
     \expandafter\@ifdefinable\csname \DMSS@prefix\DMSS@name\endcsname{%
       \def\DMSS@overwrite{00}
     }%
   \fi%
   \if\DMSS@overwrite 
     \expandafter\@firstofone%
   \else\expandafter\@gobble\fi%
   {\protected%
     \expandafter%
     \xdef\csname \DMSS@prefix\DMSS@name \endcsname{%
       \unexpanded\expandafter{\DMSS@format{#2}}%
     }}%
   \endgroup}
\define@key{DMSS}{format}{\let\DMSS@format#1}
\define@key{DMSS}{format*}{\def\DMSS@format{\expandafter#1\@firstofone}}
\define@key{DMSS}{name}{\def\DMSS@name{#1}}
\define@key{DMSS}{prefix}{\def\DMSS@prefix{#1}}
\define@key{DMSS}{overwrite}[true]{%
   \edef\DMSS@overwrite{\csname if#1\endcsname 00\else 01\fi}}
\setkeys{DMSS}{overwrite=false}
\setkeys{DMSS}{format*=}

\newcommand\MakeDeclareShorthandCommand[3]{%
   \expandafter\newcommand\expandafter*\csname math#2format\endcsname{#3}%
   \toks@{\newcommand*#1[2][]}%
   \@temptokena{,prefix=#2,##1}%
   \edef\tmpi{%
     \the\toks@{\noexpand\DeclareMathSymbolShorthand{%
     format=\expandafter\noexpand \csname math#2format\endcsname
     \the\@temptokena}{####2}}
   }%
   \tmpi
}

\MakeDeclareShorthandCommand{\DeclareMathSet}{set}{\mathbb}
\makeatother

\DeclareMathSet{R}
\DeclareMathSet{C}
\DeclareMathSet{Z}
\DeclareMathSet{Q}
\DeclareMathSet{F}
\DeclareMathSet{N}
\DeclareMathSet[prefix=group,format=\mathrm]{Diff}
\DeclareMathSet[name=Disc,format=\mathrm]{D}
\DeclareMathSet[prefix=,format=\mathrm]{Hom}
\DeclareMathSet[format=\mathrm]{Th}
\DeclareMathSet[prefix=sheaf,format=\mathcal,name=HoloMap]{O}
\DeclareMathSet[prefix=,format=\mathrm]{Princ}
\newcommand\mathbbrm[2]{\mathbb{#1}\mathrm{#2}}
\DeclareMathSet[prefix=space,format*=\mathbbrm]{CP} 

\makeatletter
\newcommand\spaceJ[1][R]{%
  J_{\@nameuse{set#1}}
}
\makeatother

\DeclareMathOperator{\ch}{ch}
\DeclareMathOperator{\Denom}{Denom}
\DeclareMathOperator{\Num}{Num}
\DeclareMathOperator{\Index}{index}
\DeclareMathOperator{\pr}{pr}
\DeclareMathOperator{\id}{id}
\DeclareMathOperator{\kernel}{ker}
\DeclareMathOperator{\Kernel}{Ker}
\DeclareMathOperator{\cokernel}{cok}
\DeclareMathOperator*{\colim}{colim}
\DeclareMathOperator{\Prim}{Prim}
\DeclareMathOperator{\proj}{proj}
\DeclareMathOperator{\Image}{Im}
\DeclareMathOperator{\Spin}{Spin}

\DeclarePairedDelimiter{\innerp}{\langle}{\rangle}

\newcommand\diff[3][\partial]{\frac{#1 #2}{#1 #3}}

\newcommand\Hodge{\mathcal{H}}

\newcommand\StdClass{\kappa}
\newcommand\bStdClass{\bar{\StdClass}}
\newcommand\bpartial{\bar{\partial}}
\newcommand\Tpi[1][\displaystyle]{{#1 T^\pi}}
\newcommand\Tlpi{ T_\pi}
\newcommand\Npi{{\displaystyle N^\pi}}
\newcommand\PontryaginPpowerOp{\mathcal{P}}

\newcommand\subprime[1][p]{_{(#1)}}
\newcommand\subprimeX[2][p]{_{(#1)#2}}
\newcommand\superprime[1][p]{^{\smash{(#1)}\vphantom{#1}}}

\numberwithin{equation}{section}

\usepackage[amsmath,thmmarks]{ntheorem}

\theoremseparator{.}
\theoremstyle{nonumberplain}
\newtheorem{conjecture}{Conjecture}
\theorembodyfont{\normalfont}

\newcommand\QedSymbol{\ensuremath{\enspace\Box}}
\newcommand\qedhere{\hfill\QedSymbol}

\theoremsymbol{\QedSymbol}
\newtheorem{proof}{Proof}
\theoremstyle{empty}
\newtheorem{proofof}{test}

\theoremstyle{plain}
\theoremsymbol{}
\theorembodyfont{\normalfont}
\newtheorem{remark*}{Remark}

\theoremstyle{plain}
\theoremnumbering{Alph}
\theorembodyfont{\itshape}
\newtheorem{theoremA}{Theorem}

\theoremnumbering{arabic}
\newtheorem{lemma}[equation]{Lemma}
\newtheorem{proposition}[equation]{Proposition}
\newtheorem{keylemma}[equation]{Key Lemma}
\newtheorem{theorem}[equation]{Theorem}
\newtheorem{corollary}[equation]{Corollary}

\usepackage{enumitem}
\setenumerate{label=\textup{(\roman*)}}

\begin{document}

\title{An integral Riemann-Roch theorem for surface bundles}
\author{Ib Madsen}
\date{January 27, 2009}
\maketitle

\section{Introduction}
\label{sec:introduction}

This paper is a response to a conjecture by T.~Akita about an integral
Riemann-Roch theorem for surface bundles, \cite{MR1892095},
\cite{MR2405898}.

Let $\pi\colon E\to B$ be an oriented surface bundle with closed
fibers of genus $g$, and equipped with a fiberwise metric. Let
$\Hodge(E)$ be the associated Hodge bundle with fibers
$H^1(E_b;\setR)$. It becomes a $g$\nobreakdash-dimensional complex
vector bundle with the complex structure induced from the Hodge star
operator, so is classified by a map from $B$ to $BU(g)$. The element
\begin{equation*}
  n!\ch_n \in H^{2n}\bigl(BU(g);\setZ\bigr)
\end{equation*}
defines a characteristic class of the Hodge bundle and we set
\begin{equation}
  \label{eq:1}
  s_n(E):= n!\ch_n\bigl(\Hodge(E)\bigr) \in H^{2n}(B;\setZ).
\end{equation}
They are torsion classes for even $n$ by \cite{MR0387496}, but not for
odd $n$ where we shall compare them with the standard
Miller-Morita-Mumford classes
\begin{equation}
  \label{eq:2}
  \StdClass_n(E):= \pi_*\bigl(c_1(\Tpi E)^{n+1}\bigr) \in H^{2n}(B;\setZ).
\end{equation}

The Grothendieck Riemann-Roch theorem, or the family index theorem of
Atiyah and Singer, yields the relation 
\begin{equation}
  \label{eq:3}
  s_{2n-1}(E)=(-1)^{n-1}(B_n/2n) \StdClass_{2n-1}(E) \quad \text{in $H^*(B;\setQ)$},
\end{equation}
cf.\ \cite{MR914849}, \cite{MR717614}. The Bernoulli numbers are
defined by the power series
\begin{equation*}
  \frac{z}{e^z-1} +\frac{z}{2} = 1+ \sum_{n=1}^\infty (-1)^{n-1}
  B_n/(2n)! z^{2n}.
\end{equation*}
Clearing denominators in \eqref{eq:3}, T.~Akita made the following
conjecture in \cite{MR1892095}:
\begin{conjecture}[Akita]
  \label{conj:akita}
  In $H^*(B;\setZ)$,
  \begin{equation*}
    \Denom(B/2n) s_{2n-1}(E) = (-1)^{n-1}\Num(B_n/2n)\StdClass_{2n-1}(E).
  \end{equation*}
\end{conjecture}
The conjecture was verified in \cite{MR2405898} for some special
values of $n$, but it turns out to be incorrect in general. A
counterexample is given in sect.~\ref{sec:disc-akit-conj} below (for
$n=p$, an odd prime), but there are many other counterexamples and the
conjecture is incorrect even in $H^*(B;\setF_p)$.

On the positive side however, one can replace the standard classes
$\StdClass_n(E)$ by another set of integral characteristic classes which
I shall denote $\bStdClass_n(E)$. In rational cohomology they differ from
$\StdClass_n(E)$ only by a sign:
\begin{equation}
  \label{eq:4}
  \bStdClass_n(E) = (-1)^n \StdClass_n(E)\quad \text{in $H^*(B;\setQ)$,}
\end{equation}
and moreover, Akita's conjecture becomes correct when we substitute
the $\StdClass_n$ classes with $\bStdClass_n$. More precisely we shall prove
\begin{theoremA}\label{thm:A}
  \begin{enumerate}
  \item In $H^{4n-2}(B;\setZ)$,
    \begin{equation*}
      2\Denom(B_n/2n) s_{2n-1}(E)= 2(-1)^n\Num(B_n/2n)\bStdClass_{2n-1}(E).
    \end{equation*}
  \item In cohomology with $p$\nobreakdash-local coefficients, the
    difference $\StdClass_n(E)-(-1)^n\bStdClass_n(E)$ is a torsion class of
    order $p$.
  \end{enumerate}
\end{theoremA}
\begin{remark*}
  The extra factor $2$ is unfortunate, and could possibly be removed
  with a more detailed consideration, cf.\ \cite{MR2219303}.
\end{remark*}

The classes $\bStdClass_n(E)$ are not as simple to define as the standard
classes $\StdClass_n(E)$. I owe the following description to Johannes Ebert,
\cite{jebert06}. Given an oriented surface bundle (with compact
fiber), and equipped with a fiberwise Riemannian metric, one has the
fiberwise  $\bpartial$\nobreakdash-operator
\begin{equation*}
  \bpartial\colon C^\infty(E;\setC)\to C^\infty(E;\Tlpi^{0,1}E).
\end{equation*}
The target consists of sections in the conjugate dual of the fiberwise
tangent bundle, i.e.\
\begin{equation*}
  \Tlpi E^{0,1}=\Hom_\setC\bigl(\overline{\Tpi E},\setC\bigr).
\end{equation*}
The index bundle of $\bpartial$ is $\setC-\Hodge(E)$. More generally,
one may twist $\bpartial$ with a complex vector bundle $W$ on $E$ to
get
\begin{equation*}
  \bpartial_W\colon C^\infty(E;W) \to C^\infty(E;\Tlpi^{0,1}E\otimes W).
\end{equation*}
The classes $\bStdClass_n(E)$ are then given by
\begin{equation}
  \label{eq:5}
  \bStdClass_n(E) =n!\ch_n(\Index \bpartial_W),
\end{equation}
where $W= \overline{\Tpi E} - \setC$ and
$\Index(\bpartial _W)$ is the (analytic) index bundle over $B$. (Its
class in $K(B)$ is independent of choice of connection on $W$ when $W$
is not holomorphic.) The index theorem gives a different description
of $\Index(\bpartial_W)$, better suited for our purpose.

Surface bundles with connected fibers of genus $g$ are classified by
$B\groupDiff(F_g)\simeq B\Gamma_g$, where $\Gamma_g$ is the mapping class
group. The proof of Theorem~\ref{thm:A} depends on the homological
description of the stable mapping class group in terms of the
cobordism space $\Omega^\infty MT(2)=\Omega^\infty\spaceCP_{-1}^\infty$
from \cite{MR2335797}. In particular the Madsen-Weiss theorem is
needed for the counterexample to Akita's conjecture.

An oriented surface bundle $\pi\colon E\to B$ induces a ``classifying map''
\begin{equation*}
  \alpha_E\colon B\to \Omega^\infty MT(2)
\end{equation*}
and the classes $\StdClass_n(E)$, $\bStdClass_n(E)$ and $s_n(E)$ are
pull\nobreakdash-backs of universal classes $\StdClass_n$, $\bStdClass_n$
and $s_n$ in the integral cohomology of $\Omega^\infty MT(2)$.  The
analysis of the relationship between these universal classes uses
infinite loop space theory and the theory of homology operations
associated with this theory.

The paper is divided up into the following sections:
\makeatletter
\setcounter{tocdepth}{1}
\renewcommand*\l@section[2]{%
  \renewcommand\numberline[1]{\item[##1.]}
  #1 
}
\begin{itemize}[label=\arabic*.,noitemsep,leftmargin=*,labelindent=\parindent,labelsep=1em]
\@starttoc{toc}%
\item[] 
\end{itemize}
\makeatother

\section{Embedded surface bundles}
\label{sec:embedd-surf-bundl}

By Whitneys's embedding theorem, a smooth manifold bundle $\pi\colon
E\to B$ with closed fibers admits a fiberwise  embedding
\begin{equation*}
  \xymatrix@C=5mm{
    \enspace E\enspace \ar@{^{(}->}[rr]^-{i} \ar[rd]^\pi& &  { B\times\setR^N } \ar[ld]^{\pr_1} 
    \\
    & B &
  }
\end{equation*}
for some $N$, and the embedding is unique up to fiberwise isotopy if
$N$ is sufficiently large.

There are obvious advantages of considering embedded fiber bundles:
pull\nobreakdash-backs become functional, and the fiberwise (or
vertical) tangent bundle $\Tpi E$ inherits an inner
product from the standard metric on euclidian space. Assuming further
that $\pi\colon E\to B$ is an \emph{oriented} surface bundle, in the
sense that $\Tpi E$ is oriented, one gets a complex structure $J\colon
\Tpi E \to \Tpi E$ by rotating vectors by $+\frac{\pi}{2}$.

The fiberwise normal bundle of the embedding,
\begin{equation*}
  \Npi E = \{ (z,v) \in E\times \setR^N \mid v \perp \Tpi E_z \},
\end{equation*}
maps to $B\times \setR^N$ by sending $(z,v)$ to
$i(z)+\bigl(\pi(z),v\bigr)$. By the tubular neighborhood theorem, we
may assume that this map restricts to an embedding of the unit disk
bundle $\setDisc(\displaystyle \Npi E)$, cf.\
\cite[\S3.2]{MR2219303}.

Given an oriented embedded surface bundle, one gets a diagram of
vector bundles
\begin{equation}
  \label{eq:6}
  \cxymatrix{
    \Npi E \ar[r]^-{\hat{t}} \ar[d]^{p} & U_{2,N-2}^\perp \ar[d] 
    \\
    E \ar[r]^-{t} & G(2,N-2) \rlap{\ ,}
  }
\end{equation}
where $G(2,N-2)$ is the Grassmann manifold of oriented
$2$\nobreakdash-planes in $\setR^N$,
\begin{align*}
  U_{2,N-2}^\perp & = \{(V,v)\mid V\in G(2,N-2), \enspace v\perp V\},
  \\
  t(z) &= TE_{\pi(z)},
  \\
  \hat{t}(z,v) &= (TE_{\pi(z)},v).
\end{align*}

Taking $N=2n+2$, we have the embedding of complex projective
$n$\nobreakdash-space $\spaceCP^n$ into $G(2,2n)$ by a
$(2n-1)$\nobreakdash-connected map, and a pull\nobreakdash-back
diagram of vector bundles
\begin{equation*}
  \xymatrix{
    L_n^\perp \ar[r] \ar[d] & U_{2,2n}^\perp \ar[d]
    \\
    \spaceCP^n \ar[r] & G(2,2n) \rlap{\ .}
  }
\end{equation*}
For $n$ large enough, $t\colon E\to G(2,2n)$ is homotopic to a map
that factors over $\spaceCP^n$, and \eqref{eq:6} can be replaced by the diagram
\begin{equation}
  \label{eq:7}
  \cxymatrix{
    \Npi E \ar[r]^-{\hat{t}} \ar[d]^{p} & L_n^\perp  \ar[d]
    \\
    E \ar[r]^-{t} &  \spaceCP^n  \rlap{\ .}
  }
\end{equation}
One consequence of \eqref{eq:7} is that we may assume that
$\displaystyle \Npi E$ is a complex vector bundle. Alternatively a
complex structure can be seen as follows: if $E\subset
B\times\setR^{n+2}$ is a fiberwise embedding with vertical normal
bundle $\displaystyle \Npi E $, then $\Tpi E \subset
B\times T\setR^{n+2}$ has complex normal bundle $q^*(\displaystyle
\Npi E \otimes_{\setR} \setC)$, where $q$ is the bundle projection of
$\Tpi E$, \cite{MR0236950}, \cite{MR0279833}, and $E\subset \Tpi E
\subset B\times T\setR^{n+2}$ has normal bundle $\Tpi E \oplus (\Npi
E) \otimes_{\setR} \setC$ ; this is a complex vector bundle.

For the rest of the paper, a surface bundle $\pi\colon E\to B$ will
mean an embedded oriented surface bundle $E\subset B\times
\setR^{2n+2}$ with compact fiber and complex vertical normal bundle,
such that the standard map $\Npi E\to B\times \setR^{2n+2}$ restricts to
an embedding of the unit disc bundles $\setDisc (\Npi E)$. The integer
$n$ is not part of the structure: $E\subset B\times \setR^{2n+2}$ is
identified with $E\subset B\times \setR^{2n+4}$.

The diffeomorphism classes of such surface bundles are in bijective
correspondence with the set of homotopy classes $[B,B\groupDiff(F)]$, where
$F$ is the fiber surface and $\groupDiff(F)$ denotes the topological group
of orientation preserving diffeomorphisms, cf.\ \S2.1 of \cite{art18}.

For such a surface bundle, we can apply the Pontryagin-Thom
construction on the embedding
\begin{equation*}
  \setDisc(\Npi E) \hookrightarrow B\times \setR^{2n+2},
\end{equation*}
giving the map
\begin{equation}
  \label{eq:8}
  c_E\colon B_+\wedge S^{2n+2} \to \setTh(\Npi E) := \setDisc (\Npi
  E)/\partial\ . 
\end{equation}
There is an induced push\nobreakdash-forward map $\pi_*$ in
$K$\nobreakdash-theory defined by the diagram
\begin{equation*}
  \xymatrix{
    K(E) \ar[r]^-{\Phi} \ar[d]^{\pi_*} 
    & 
    { \widetilde{K}\bigl(\setTh(\Npi
    E)\bigr) } 
    \ar[d]^{(c_E)^*}
    \\
    K(B) \ar[r]^-{\Phi}_-{\cong} 
    & 
    { \widetilde{K}(B_+\wedge S^{2n+2})
      }\rlap{\ ,}
  }
\end{equation*}
where $\Phi$ denotes the $K$\nobreakdash-theory Thom isomorphism. It
depends on the choice of Thom class; we use the class from
\cite[\S2.7]{MR1043170}. 

Let $\ch_n\colon K(B) \to H^{2n}(B;\setQ)$ be the $n$'th component of
the Chern character, and
\begin{equation*}
  s_n = n!\ch_n\colon K(B) \to H^{2n}(B;\setZ)
\end{equation*}
the accompanying map into integral cohomology. The classes that appear
in Theorem~\ref{thm:A} have the following alternative description
\begin{equation}
  \label{eq:9}
  \begin{split}
    s_n(E) &= (-1)^n s_n\bigl(\pi_*(1)\bigr),
    \\
    \bStdClass_n(E) &= s_n\bigl(\pi_*([\overline{\Tpi E}]-1)\bigr),
  \end{split}
\end{equation}
where $\overline{\Tpi E}$ is conjugate to $(\Tpi E, J)$ and isomorphic
to the vertical complex cotangent bundle $\Hom_{\setC} (\Tpi
E,\setC)$.

The relationship between \eqref{eq:9} and the description given in the
introduction is the index theorem \cite{MR0279833}. This is discussed
in the next section.

\section{The fiberwise $\bm{\bpartial}$\protect\nobreakdash-operator and its index}
\label{sec:fiberw-bpart-oper}

This short section explains the family index theorem for the fiberwise
$\bpartial$\nobreakdash-operator. This is well\nobreakdash-known for
the experts, but I will include some details for the more
inexperienced reader.

Let $V$ be a vector space with inner product and isometric complex
structure $J\colon V\to V$, $J^2=-\id$. The involution $\phi\mapsto
i\phi J^{-1}$ on $\Hom_{\setC}(V,\setC)$ decomposes into $\pm1$
eigenspaces
\begin{equation*}
  \Hom_{\setR}(V,\setC) = \Hom_{\setC}(V,\setC) \oplus \Hom_{\setC}(\overline{V},\setC)
\end{equation*}
where on the right hand side $V=(V,J)$ and $\overline{V} =(V,J^{-1})$.
The subspace $\Hom_{\setR}(V,\setR)$ sits diagonally in the direct
sum, and the projections define isomorphisms
\begin{equation}
  \label{eq:10}
  \Hom_{\setR}(V,\setR)\cong \Hom_{\setC}(V,\setC), \quad 
  \Hom_{\setR}(V,\setR)\cong \Hom_{\setC}(\overline{V},\setC).
\end{equation}
The two isomorphisms are given by $\phi\mapsto 1/2(\phi+i\phi J^{-1})$
and $\phi\mapsto 1/2(\phi-i\phi J^{-1})$, respectively. The inner
product $\innerp{\,,\,}$ on $V$ identifies $V$ with
$\Hom_{\setR}(V,\setR)$ and the combined isomorphism
\begin{equation}
  \label{eq:11}
  V\xrightarrow{\ \cong\ } \Hom_{\setR}(V,\setR)
  \xrightarrow{\ \cong\ } \Hom_{\setC}(\overline{V},\setC),
\end{equation}
which sends $w\in V$ to $1/2(\phi_w-i\phi_wJ^{-1})$ with
$\phi_w(v)=\innerp{v,w}$, is a complex isomorphism.

Applied fiberwise to the vertical tangent bundle $\Tpi E$ of an
embedded surface bundle, define
\begin{equation*}
  T_\pi ^{1,0} E := \Hom_{\setC}(\Tpi E, \setC), \quad 
  T_\pi^{0,1} E := \Hom_{\setC}(\overline{\Tpi E}, \setC)
\end{equation*}
so that
\begin{equation*}
  C^\infty(E,\Tpi E) = C^\infty(E,T_\pi^{1,0} E)
    \oplus
    C^\infty(E,T_\pi^{0,1} E).
\end{equation*}

The fiberwise exterior differential from $C^\infty(E,\setC)$ into
$C^\infty(E,\Tpi E)$ decomposes accordingly. The second component is
\begin{equation*}
  \bpartial\colon C^\infty(E,\setC) \to C^\infty(E,T_\pi^{0,1} E),
\end{equation*}
which on a fiber $F$ of $\pi\colon E\to B$ is given locally by
$\bpartial\phi=\diff{\phi}{\bar{z}} d\bar{z}$ with respect to a
complex coordinate on $F$ compatible with $J$.

More generally, if $W$ is a complex vector bundle on $E$ then
(possibly after choice of connection) one has the associated operator
\begin{equation*}
  \bpartial_W\colon C^\infty(E,W) \to C^\infty(E,T_\pi^{0,1} E \otimes
  W). 
\end{equation*}
The symbol $\sigma(\bpartial_W)$ is identified as
\begin{equation}
  \label{eq:12}
  \sigma(\bpartial_W) = \Phi([W])
  =  \lambda_{\Tpi[]\! E} \cdot [W] \in \tilde{K}\bigl(\setTh(\Tpi E)\bigr),
\end{equation}
cf.\ \cite[\S4]{MR0236952}. This class is independent of choice of
connection. The (virtual) index bundle is
\begin{equation*}
  \Index(\bpartial_W) =\kernel(\bpartial_W) - \cokernel(\bpartial_W).
\end{equation*}
It defines a class in $K(B)$, and the index theorem of \cite{MR0279833}
asserts that
\begin{equation}
  \label{eq:13}
  [\Index(\bpartial_W)] = \pi_!\bigl(\sigma(\bpartial_W)\bigr) \in K(B).
\end{equation}
The push\nobreakdash-forward $\pi_!\colon \tilde{K}\bigl(\setTh(\Tpi
E)\bigr) \to K(B)$ of \cite{MR0279833} is related to our $\pi_*\colon
K(E) \to K(B)$ via the formula
\begin{equation}
  \label{eq:14}
  \pi_*([W]) =\pi_!\bigl(\Phi([W])\bigr).
\end{equation}
Hence \eqref{eq:13} has the alternative formulation:
\begin{equation}
  \label{eq:15}
  [\Index\bpartial_W] = \pi_*([W]),
\end{equation}
and the classes of \eqref{eq:9} are given by
\begin{equation}
  \label{eq:16}
  s_n(\pi_*[W]) = s_n(\Index \bpartial_W),
\end{equation}
for $[W]=1$ and $[W] = [\overline{\Tpi E}] -1$, respectively. We have
left to identify $\Index(\bpartial)$ in terms of the Hodge bundle.
The result is:
\begin{lemma}\label{lem:3.4}
  \begin{enumerate}
  \item \label{item:1} $\kernel \bpartial = B\times \setC$,
  \item \label{item:2} $\cokernel \bpartial =\Hodge (E)$, the Hodge
    bundle, and $\pi_*(1) = 1-[\Hodge(E)]$.
  \end{enumerate}
\end{lemma}
\begin{proof}
  On a fiber $F$ of $\pi\colon E\to B$, the kernel of $\bpartial_{|F}$
  consists of the holomorphic maps from $F$ to $\setC$ which, since
  $F$ is a closed Riemannian surface, are the constant maps. This
  proves \ref{item:1}.  The cokernel of $\bpartial_{|F}$ is
  \begin{equation*}
    \cokernel \bpartial_{|F} = H^1(F;\sheafHoloMap),
  \end{equation*}
  where $\sheafHoloMap$ is the sheaf of holomorphic maps, see e.g.\
  \cite[Thm.~15.14]{MR648106}. On the other hand, by Serre duality
  \begin{equation*}
    H^1(F;\sheafHoloMap) =\Hom_{\setC}\bigl(\Omega(F),\setC\bigr)
  \end{equation*}
  where $\Omega(F)$ is the vector space of holomorphic
  $1$\nobreakdash-forms, \cite[Thm.~17.9]{MR648106}.

  We used the Hodge $*$\nobreakdash-operator to give the space
  $H^1(F;\setR)$, of harmonic $1$\nobreakdash-forms, a complex
  structure, isometric w.r.t.\ the inner product given by Poincar\'e
  duality. The isomorphism
  \begin{equation*}
    H^1(F;\setR) \xrightarrow{\ \cong\ } \Omega(F),
  \end{equation*}
  given locally by sending the harmonic $1$\nobreakdash-form $\phi
  dx+\psi dy$ to $(\phi -i\psi)dz$, is conjugate linear. This yields
  the complex linear isomorphism
  \begin{equation*}
    \bigl(H^1(F;\setR),*\bigr) \cong \Hom_{\setC} \bigl(\Omega(F),\setC\bigr).
  \end{equation*}
  Since the Hodge bundle is
  \begin{equation*}
    \Princ(E)\times_{\groupDiff(F)} \bigl(H^1(F;\setR),*\bigr)
  \end{equation*}
  this proves \ref{item:2}.
\end{proof}

\section{The main diagram}
\label{sec:main-diagram}

Let $E\subset B \times \setR^{2n+1}$ be an embedded, oriented surface
bundle as in \S\ref{sec:introduction}, classified by
\begin{equation*}
  \xymatrix{
    { \Npi E } \ar[r]^-{\hat{t}} \ar[d] & { L_n^\perp } \ar[d]
    \\
    E \ar[r]^-{t} & \spaceCP^n  \rlap{\ ,}
  }
\end{equation*}
and inducing the map
\begin{equation*}
  \hat{t} \colon \setTh(\Npi E) \to \setTh(L_n^\perp).
\end{equation*}
Let
\begin{equation}
  \label{eq:17}
  \omega\colon \setTh(L_n^\perp) \to \setTh(L_n \oplus L_n^\perp) =
  \spaceCP_+^n \wedge S^{2n+2}
\end{equation}
be the inclusion along the zero section of $L_n$. The
elements $\pi_*(1)$ and $\pi_*(\overline{\Tpi E} -1)$ of $K(B)$ are
explicitly given by the formulas
\begin{equation}
  \label{eq:18}
  \begin{split}
    \pi_*(1) &= \Phi^{-1}_B c_E^* \hat{t}^*(\lambda_{L_n^\perp}),
    \\
    \pi_*(\overline{\Tpi E}-1) &=
    \Phi^{-1}_Bc_E^*\hat{t}^*\omega^*\Phi_{\spaceCP^n}(\bar{L}_n),
  \end{split}
\end{equation}
where $\Phi_X\colon K(X)\to \tilde{K}(X_+\wedge S^{2n+2})$ is the Thom
isomorphism for the trivial bundle on $X$, and $c_E\colon B_+\wedge
S^{2n+2} \to \setTh(\Npi E)$ is the collapse map of \eqref{eq:8}.
For the second formula in \eqref{eq:18}, note that
$\omega^*(\lambda_{L+L^\perp}) = (1-L)\lambda_L^\perp$.

We now describe the key relation between the elements
$\lambda_{L_n^\perp}$and $\omega^*\Phi_{\spaceCP^n}(\bar{L}_n)$ which
eventually leads to the proof of Theorem~\ref{thm:A}. 

Given an $n$\nobreakdash-dimensional complex vector bundle $W$ on $X$,
\begin{equation*}
\varrho^k(W)\in K(X)\otimes \setZ[\tfrac{1}{k}]
\end{equation*}
is the $K$\nobreakdash-theory characteristic class defined by
\begin{equation*}
  \psi^k(\lambda_W)=k^n\varrho^k(W)\lambda_W,
\end{equation*}
where $\psi^k$ is the $k$'th Adams operation. Then
\begin{align*}
  \varrho^k(V\oplus W) &= \varrho^k(V) \varrho^k(W),
  \\
  \varrho^k(L) &= \tfrac{1}{k}(1+L+\dots+L^{k-1}), \quad \text{$L$ a
    line bundle.}
\end{align*}
Define the operation $r^k$ by the diagram
\begin{equation}
  \label{eq:19}
  \cxymatrix[@C=1.5cm]{
    K(X) \ar[r]^-{\Phi}_-{\cong} \ar[d]^{r^k} 
    & \tilde{K}(X_+\wedge S^2) \ar[d]^{\varrho^k}
    \\
    K(X)\otimes \setZ[\tfrac{1}{k}] \ar[r]^-{1+\Phi}_-{\cong} &
    1+\tilde{K}(X_+\wedge S^2)\otimes \setZ[\tfrac{1}{k}]
    \rlap{\ .}
  }
\end{equation}
It is additive,
\begin{equation*}
  r^k(V\oplus W) =r^k(V) + r^k(W),
\end{equation*}
and for a line bundle $L$ it is given by the formula (with $\bar{L}$
the conjugate line bundle)
\begin{equation}
  \label{eq:20}
  r^k(\bar{L}) = \frac{L^{k-2}+ 2L^{k-3}+\dots+(k-2)L+(k-1)}%
                      {1+L+\dots+L^{k-1}} \ .
\end{equation}
This follows easily upon using that $\Phi$ in \eqref{eq:18} is
(exterior) multiplication by $\lambda_{\setC}=L_1-1$ in
$\tilde{K}(S^2)$, together with the equation
\begin{equation*}
  (L-1)\bigl(L^{k-2}+2L^{k-3} +\dots+ (k-1)\bigr)
  =(1+L+\dots+L^{k-1})-k,
\end{equation*}
cf.\ \cite[Prop.~6.2]{MR1778114}. Since $\omega^*(\lambda_{L_n \oplus
  L_n^\perp})=(1-L_n)\lambda_{L_n^\perp}$, we get
\begin{equation}
  \label{eq:21}
  (k^{-n}\psi^k-\id)(\lambda_{L_n^\perp})
  =\omega^*\bigl(r^k(\bar{L}_n)\lambda_{L_n \oplus L_n^\perp}\bigr).
\end{equation}
This relation can also be expressed as the homotopy commutative
diagram below, where $(\setZ\times\nobreak BU)[\tfrac{1}{k}]$ is the
classifying space of $K(X)\otimes \setZ[\tfrac{1}{k}]$, i.e.\
\begin{equation*}
  K(X)\otimes \setZ[\tfrac{1}{k}] = \bigl[X,(\setZ \times
  BU)[\tfrac{1}{k}]  \bigr]
\end{equation*}
for compact $X$. The diagram is
\begin{equation}
  \label{eq:22}
  \cxymatrix[@C=2cm]{
    \setTh(L_n^\perp) \ar[r]^-{\omega} \ar[d]^{\lambda_{L_n^\perp}}
    &
    \spaceCP_+^n \wedge S^{2n+2} \ar[r]^-{\bar{L}_n \wedge S^{2n+2}}
    &
    (\setZ \times BU)\wedge S^{2n+2} \ar[d]^{r^k \wedge S^{2n+2}}
     \\
    { \setZ\times BU } \ar[r]^-{k^{-n}\psi^k-\id}
    &
   { (\setZ\times BU)[\tfrac{1}{k}] }
    &
   { (\setZ \times BU)[\tfrac{1}{k}] \wedge S^{2n+2} }
      \ar[l]_-{\varepsilon} \rlap{\ ,}
  }
\end{equation}
with $\varepsilon\colon (\setZ\times BU) \wedge S^{2n+2} \to BU$
induced from multiplication with $\lambda_{\setC^{n+1}} \in
\tilde{K}(S^{2n+2})$.

In order to make the above diagram independent of the embedding
dimension, we pass to spectra and their associated infinite loop
spaces, cf.\ \cite{MR1324104}, \cite{MR505692}.

We remember that for a (pre) spectrum $A=\{A_n,\varepsilon_n\}$ the
associated infinite loop space is
\begin{equation*}
  \Omega^\infty A = \colim \Omega^n A_n
\end{equation*}
where the maps $\Omega^nA_n \to \Omega^{n+1}A_{n+1}$ comes from the
adjoint $\varepsilon_n'\colon A_n \to \Omega A_{n+1}$.  If $A$ is an
($\Omega$-) spectrum, i.e.\ if $\varepsilon_n'\colon A_n\to\Omega
A_{n+1}$ is a homotopy equivalence, then $\Omega^\infty A \simeq A_0$.
The spectra relevant to us are $\Sigma^\infty Y$, $KU$ and $MT(2)$.
They have $(2n+2)$'nd spaces:
\begin{align*}
  (\Sigma^\infty Y)_{2n+2} &= Y\wedge S^{2n+2}, 
  \\
  KU_{2n+2} &= \setZ \times BU,
  \\
  MT(2)_{2n+2} &= \setTh(L_n^\perp).
\end{align*}
Applying $\Omega^{2n+2}(-)$ to \eqref{eq:22} and taking colimit over
$n$ leads to the following homotopy commutative diagram of infinite
loop spaces, where we use the standard notations: $Q(Y)=\Omega^\infty
\Sigma^\infty Y$ and $X_+=X\sqcup \{+\}$,
\begin{equation}
  \label{eq:23}
  \cxymatrix[@C=15mm]{
    \Omega^\infty MT(2) \ar[r]^-{\omega} \ar[d]^{\lambda_{-L}}
    &
    Q(\spaceCP^\infty_+) \ar[r]^-{Q(\bar{L})} 
    &
    Q(\setZ\times BU) \ar[d]^{Q(r^k)}
    \\
    \setZ \times BU \ar[r]^-{k\psi^k-\id} 
    &
    \setZ \times BU[\tfrac{1}{k}] 
    & 
    Q\bigl((\setZ \times BU)[\tfrac{1}{k}]\bigr)
    \ar[l]_-{E} \rlap{\ .}
  }
\end{equation}
The map $E\colon Q(\setZ \times BU) \to \setZ\times BU$, and its
localised version with $k$ inverted, is a consequence of Bott
periodicity $\Omega^{2n}(\setZ \times BU)\simeq \setZ \times BU$. More
generally, for any $\Omega$\nobreakdash-spectrum
$A=\{A_n,\varepsilon_n\}$, the maps
\begin{equation*}
  \Omega^n(A_0\wedge S^n) \to \Omega^n\! A_n \xleftarrow{\:\simeq\:} A_0 
\end{equation*}
determines a (weak) homotopy class $E\colon Q(A_0)\to A_0$.

\penalty -500 

The composition
\begin{equation*}
\hat{t} \circ  c_E\colon B_+\wedge S^{2n+2} \to \setTh(\Npi E) \to
  \setTh(L_n^\perp) 
\end{equation*}
of an embedded surface bundle is the $(2n+2)$'nd component of a map of
spectra
\begin{equation*}
  c_E\colon \Sigma^\infty(B_+)\to MT(2).
\end{equation*}

Its adjoint is
\begin{equation}
  \label{eq:24}
  \alpha_E\colon B \xrightarrow{\: i_B\:} Q(B_+)
  \xrightarrow{\:\Omega^\infty c_E\:} \Omega^\infty MT(2),
\end{equation}
where $i_B$ is the inclusion.

With these notations, the two elements of \eqref{eq:18} are represented
by the homotopy classes
\begin{equation}
  \label{eq:25}
  \begin{split}
    &B\xrightarrow{\:\alpha_E\:} \Omega^\infty MT(2) 
    \xrightarrow{\: \lambda_{-L} \:} \setZ\times BU,
    \\
    & B\xrightarrow{\:\alpha_E\:} \Omega^\infty MT(2) 
    \xrightarrow{\: \omega \:} 
    Q(\spaceCP_+^\infty) 
    \xrightarrow{\:Q(\bar{L})\:}
    Q(\setZ\times BU) \xrightarrow{\: E \:} \setZ \times BU.
  \end{split}
\end{equation}

Theorem~\ref{thm:A} is proved by evaluating \eqref{eq:23} on
cohomology, but to accomplish this one needs to examine the diagram
\begin{equation}
  \label{eq:26}
  \cxymatrix{
    Q(\setZ\times BU) \ar[r]^-{E} \ar[d]^{Q(r^k)}
    &
    \setZ\times BU \ar[d]^{r^k}
    \\
    Q\bigl((\setZ\times BU)[\tfrac{1}{k}]\bigr) 
    \ar[r]^-{E}
    &
    (\setZ\times BU)[\tfrac{1}{k}]
  }
\end{equation}
on the cohomological level. Unfortunately, \eqref{eq:26} is \emph{not}
a homotopy commutative diagram: $r^k$ is only twice deloopable but not
an infinite loop map.

We may restrict \eqref{eq:26} to the zero component $BU=\{0\}\times BU
\subset \setZ \times BU$, and localise at a prime $p$ with $(k,p)=1$.

\begin{keylemma}\label{lemma:key}
  For $(k,p)=1$, the cohomological diagram
  \begin{equation*}
    \xymatrix{
      \Prim H^*\bigl( Q(BU_{(p)});\setZ \bigr)
      & \Prim H^*(BU_{(p)};\setZ) \ar[l]_-{E^*}
      \\
      \Prim H^* \bigl(Q(BU_{(p)});\setZ\bigr) \ar[u]^{Q(r^k)^*}
      &
      \Prim H^*(BU_{(p)};\setZ)
      \ar[l]_-{E^*}
      \ar[u]^{(r^k)^*}
    }
  \end{equation*}
  is homotopy commutative.
\end{keylemma}

Here $\Prim(\ )$ denotes the primitive elements in the Hopf algebras,
e.g.\
\begin{equation*}
  \Prim H^{2n}(BU_{(p)};\setZ) =\setZ_{(p)} \innerp{s_n}, \quad
  s_n=n!\ch_n, 
\end{equation*}
where $\setZ_{(p)} \subset \setQ$ are the fractions with denominator
prime to $p$. The lemma is proved in the next section.

\section{The key lemma}
\label{sec:key-lemma}

This section deals with the non-commutativity of diagram
\eqref{eq:26}. Our method is to use the Artin-Hasse logarithm of
\cite{MR1000393}. To this end we single out a prime $p$ and only
consider $\varrho^k$ and $r^k$ for $(k,p)=1$. Let us write $K\subprime
(X) = K(X)\otimes \setZ\subprime$.  We have the group homomorphism
\begin{equation*}
  \varrho^k\colon \tilde{K}\subprime(X) \to 1+ \tilde{K}\subprime(X)
\end{equation*}
where the group structure on the target is tensor product of virtual
vector bundles of dimension~$1$.

The representing space for $\tilde{K}(X)$ is $BU^\oplus =
\{0\}\times BU$. This is the infinite loop space associated to the
connected $K$\nobreakdash-theory spectrum $bu^\oplus$ of \cite{MR1324104}. The
representing space for $1+\tilde{K}(X)$ is $BU^\otimes =\nobreak
\{1\}\times\nobreak BU$. This is the infinite loop space of a spectrum
$bu^\otimes$, constructed from the abstract theory of infinite loop
spaces \cite{MR0420609}, \cite{MR0420610}, \cite{MR0353298}. The map
\begin{equation*}
  \varrho^k\colon BY\subprime ^\oplus \to BU\subprime^\otimes
\end{equation*}
is an infinite loop map by \cite{MR0494076}, i.e.\ it lifts to a map
of spectra $bu\subprime^\oplus \to bu\subprime^\otimes$.

The Artin-Hasse logarithm
\begin{equation*}
  L\subprime\colon 1+\tilde{K}\subprime(X) \to \tilde{K}\subprime(X)
\end{equation*}
is for compact $X$ defined by the formula
\begin{equation*}
  L\subprime(1-x) = -\sum_{(n,p)=1} \frac{1}{n}\sum_{t=0}^\infty
  \theta^{p^t}(x^n), 
\end{equation*}
where $\theta^{p^t}\colon \tilde{K}\subprime(X)\to
\tilde{K}\subprime(X)$ is the operation
\begin{equation*}
  \theta^{p^t}(x) = \frac{1}{p^t} \bigl(x^{p^t}
  -\psi^p(x^{p^{t-1}})\bigr), \quad t>0
\end{equation*}
and $\theta^1(x)=x$. It exists in $\tilde{K}\subprime(X)$ because
$\psi^p$ is multiplicative and $\psi^p(x) \equiv x^p\pmod{p}$, and it
is uniquely defined since $\tilde{K}(BU)$ is torsion free. Note that
rationally,
\begin{equation*}
  L\subprime(1-x) = \left(\frac{\psi^p}{p}-1\right)\log(1-x) \in
  \tilde{K}(X)\otimes \setQ.
\end{equation*}
It is the double loop of $L\subprime$ we
are interested in, i.e.\ in the map $l\subprime$ defined by the diagram
\begin{equation*}
  \xymatrix{
    1+\tilde{K}\subprime(X\wedge S^2) \ar[r]^-{L\subprime} 
    &
    \tilde{K}\subprime(X\wedge S^2)
    \\
    \tilde{K}\subprime(X) \ar[r]^-{l\subprime}
    \ar[u]^{\cong}_{1+\Phi}
    & 
    \tilde{K}\subprime(X) \rlap{\ .}
    \ar[u]^{\cong}_{\Phi}
  }
\end{equation*}

\begin{lemma}\label{lem:5.1}
  The map $l\subprime$is given by the formula
  \begin{equation*}
    l\subprime(x) = x+\psi^p(x).
  \end{equation*}
\end{lemma}
\begin{proof}
  All products vanish in reduced $K$\nobreakdash-theory of a
  suspension, so
  \begin{align*}
    -L\subprime(1-x\lambda_{\setC}) = \smash[b]{ \sum_{t=0}^\infty }
    \psi^{p^t}(x\lambda_{\setC}) & = x\lambda_{\setC} +\tfrac{1}{p}
    \psi^p(x\lambda_{\setC})
    \\
    &= x\lambda_{\setC} + \psi^p(x)\varrho^p(1)\lambda_{\setC}
    \\
    &= \bigl(x+\psi^p(x)\bigr)\lambda_{\setC}
  \end{align*}
  Strictly speaking, this calculation only makes sense when
  $K\subprime(X)$ is torsion free. But this is the case when $X=BU$.
  In general the formula follows by naturality from this case.
\end{proof}
Based on the criteria from \cite{MR0494076}, tom Dieck showed in
\cite{MR1000393} that 
\begin{equation*}
  L\subprime\colon BSU\subprime^\otimes \to BSU\subprime^\oplus
\end{equation*}
is an infinite loop map. The composite $L\subprime \circ \varrho^k$
from $BSU\subprime^\oplus$ to itself is then also infinitely
deloopable. Taking double loops implies that
\begin{equation*}
  BU\subprime^\oplus \xrightarrow{\:r^k\:} BU\subprime
  \xrightarrow{\:l\subprime\:} BU\subprime^\oplus
\end{equation*}
is an infinite loop map, and the diagram
\begin{equation}
  \label{eq:27}
  \cxymatrix[@C=12mm]{
    Q(BU\subprime) \ar[r]^-{Q(r^k)} \ar[d]^{E}
    &
    Q(BU\subprime) \ar[r]^-{Q(l\subprime)} 
    &
    Q(BU\subprime)\ar[d]^{E}
    \\
    BU\subprime \ar[r]^-{r^k} 
    &
    BU\subprime \ar[r]^-{l\subprime} 
    & 
    BU\subprime
  }
\end{equation}
is consequently homotopy commutative. Lemma~\ref{lem:5.1} shows that
$l\subprime$ induces isomorphism on homotopy groups, so is a homotopy
equivalence.

\begin{lemma}\label{lem:5.3}
  The maps $l\subprime$ and $Q(l\subprime)$ induce the identity on
  cohomology with $\setF_p$ coefficients.
\end{lemma}
\begin{proof}
  It suffices to check that
  \begin{equation*}
    l\subprimeX*\colon H_*(BU\subprime;\setF_p) \to
    H_*(BU\subprime;\setF_p) 
  \end{equation*}
  is the identity, since $H_*\bigl(Q(X);\setF_p\bigr)$ is a functor of
  $H_*(X,\setF_p)$. 

  The Adams operation $\psi^p$ induces multiplication by $p^n$ on
  $H_{2n}(BU\subprime;\setZ)$, so induces the zero map on
  $\tilde{H}_*(BU\subprime;\setF_p)$. By Lemma~\ref{lem:5.1},
  $l\subprime=\id$.
\end{proof}

We need an integral or $p$\nobreakdash-local version of the previous
lemma. This requires some infinite loop space theory, namely the
Pontryagin $p$'th power homology operations (Frobenious operations),
\begin{equation*}
  \PontryaginPpowerOp\colon H_n(X;\setZ/p^r) \to
  H_{pn}(X;\setZ/p^{r+1}). 
\end{equation*}
They are defined when $X$ is infinite loop and are natural with
respect to maps in that category. The key property is that
$\varrho_r\PontryaginPpowerOp(x) =x^p$, where $\varrho_r$ denotes
reduction to $\setZ/p^r$ coefficients. The reader is referred to
\cite{MR1778114}, \cite{MR0436146} and the present appendix for
properties of $\PontryaginPpowerOp$.

If $X$ is a space with multiplication $m\colon X\times X\to X$, e.g.\
an infinite loop space, a cohomology class $\gamma\in
H^*(X;\setZ\subprime)$ is called \emph{primitive} if
\begin{equation*}
  m^*(\gamma)=\mu(1\otimes \gamma+\gamma\otimes1)
\end{equation*}
where
\begin{equation*}
  \mu\colon H^*(X;\setZ\subprime)\otimes H^*(X;\setZ\subprime) \to 
  H^*(X \times X;\setZ\subprime) 
\end{equation*}
is the exterior product.

\begin{theorem}\label{thm:5.5}
  Let $X=Q(Y)$ where $Y$ is a space of finite type, and let $\gamma\in
  H^*(X;\setZ\subprime)$ be a primitive cohomology class such that
  \begin{enumerate}
  \item \label{item:7} $\varrho_1(\gamma)=0$ in $H^*(X;\setZ/p)$,
  \item \label{item:8} $i^*(\gamma)=0$, where $i$ is the inclusion $i\colon Y\to
    Q(Y)$,
  \item \label{item:9} $\innerp{\gamma,\PontryaginPpowerOp\superprime[r](c)}=0$ in $\setZ/p^{r+1}$ for
    all $c\in H_*(X;\setZ/p)$.
  \end{enumerate}
  Then $\gamma=0$.\qedhere
\end{theorem}
The proof is deferred to the Appendix.

\begin{corollary}\label{cor:5.6}
  Let $l\subprime$ and $Q(l\subprime)$ be the maps from \eqref{eq:27},
  and $s_n=n!\ch_n \in H^{2n}(BU\subprime;\setZ)$. Then
  \begin{equation*}
    l\subprime^*(s_n)=(1+p^n)s_n
    \quad\text{and}\quad
    Q(l\subprime)^*E^*(s_n)=(1+p^n)E^*(s_n).
  \end{equation*}
\end{corollary}
\begin{proof}
  The first equation follow because $s_n$ is primitive and
  $\psi^p(s_n)=p^ns_n$. We have $Q(l\subprime)=\id+Q(\psi^p)$, so it
  suffices to show that
  \begin{equation*}
    \gamma=p^nE^*(s_n)-Q(\psi^p)^*E^*(s_n)
  \end{equation*}
  is the zero class in $H^{2n}\bigl(Q(BU\subprime);\setZ\bigr)$. We
  check the three conditions of Theorem~\ref{thm:5.5}. The first
  condition is satisfied because $i^*Q(l\subprime)=l\subprime$, the
  second condition by Lemma~\ref{lem:5.3}. 

  The third condition that $\gamma$ vanishes on iterated powers of
  $\PontryaginPpowerOp$ is also satisfied. Indeed, both $E$ and
  $Q(\psi^p)$ are infinite loop maps, so commute with
  $\PontryaginPpowerOp$. Moreover each $x\in H_*(BU\subprime;\setZ/p)$
  is the reduction of an integral class $\hat{x}$, so
  $\PontryaginPpowerOp\superprime[r](x)=\varrho_{r+1}(\hat{x})^{p^r}$
  which is annihilated by the primitive class~$s_n$.
\end{proof}

The Key Lemma~\ref{lemma:key} is an obvious consequence of
diagram~\eqref{eq:27} and Corollary~\ref{cor:5.6}.

\section{Proof of Theorem~\ref{thm:A}}
\label{sec:proof-theor-refthm:a}

Our first objective is to evaluate
\begin{equation*}
  (r^k)^*\colon \Prim H^*(BU\subprime;\setZ) \to \Prim H^*(BU\subprime;\setZ).
\end{equation*}
We begin with some notation from \cite{MR0198468}. Let $W$ be a
complex vector bundle on $X$ of dimension~$m$. Define $bh(W)\in
H^*(X;\setQ)$ by
\begin{equation}
  \label{eq:28}
  \ch \lambda_W = (-1)^m bh(W)U_W,
\end{equation}
where $U_W\in H^{2m}(\setTh(W);\setZ)$ is the cohomological Thom
class. If $m=1$, $bh(W)=\frac{e^x-1}{x}$, where $x$ is the first Chern
class of $W$. 

Let $\psi^k_H$ be the endomorphism of $H^{2*}(X;\setZ)$ that
multiplies by $k^n$ in degree $2n$, so that $\mathord{\ch}\circ \psi^k = \psi^k_H\circ
\ch$.
We need the relations
\begin{equation}
  \label{eq:29}
  \begin{split}
    \psi^k_H\bigl(bh(W)\bigr) &= \ch\bigl(\varrho^k(W)\bigr) bh(W)
    \\
    \log\bigl(bh(W)\bigr) &= \tfrac{1}{2} \ch_1(W)+\sum_{n=1}^\infty
    (-1)^{n-1}(B_n/2n) \ch_{2n}(W).
  \end{split}
\end{equation}
The first one follows by applying the Chern character to the defining
relation $\psi^k(\lambda_W)=k^m\varrho^k(W)\lambda_W$; the second
relation is from \cite[Lemma~2.4]{MR0198468}.

\begin{lemma}\label{lem:6.3}
  The map
  \begin{equation*}
    (r^k)^*\colon \Prim H^{4n-2}(BU\subprime;\setZ) \to \Prim
    H^{4n-2}(BU\subprime ;\setZ)
  \end{equation*}
  is multiplication by $(-1)^{n-1}(k^{2n}-1)B_n/2n$.
\end{lemma}
\begin{proof}
  The Hurewicz homomorphism onto the module of indecomposable elements
  \begin{equation*}
    \tilde{K}(S^{2m}) \to QH_{2m}(BU;\setZ)
  \end{equation*}
  is non\nobreakdash-zero, so the paring
  \begin{equation*}
    \tilde{K}(S^{2m})\otimes \Prim H^{2m}(BU;\setZ) \to \setZ
  \end{equation*}
  is non\nobreakdash-degenerate. Thus it suffices to evaluate
  $(r^k)_*$ on $\tilde{K}\subprime(S^{4n-2})$, or equivalently,
  $(\varrho^k)_*$ on $\tilde{K}\subprime(S^{4n})$.  We use that
  \begin{equation*}
    \ch_{n}\colon \tilde{K}(S^{2n}) \to H^{2n}(S^{2n};\setZ)
  \end{equation*}
  is an isomorphism. This is clear for $n=1$ and follows in general
  since $\ch$ preserves exterior multiplication.

  Let $u\in\tilde{K}\subprime(S^{4n})$ and set $\varrho^k(u) =
  \lambda_{2n}\cdot u$. The first equation in \eqref{eq:29} gives
  \begin{align*}
    \psi^K_H\bigl(\log bh(u)\bigr) -\log bh(u) & = \log \ch
    \bigl(\varrho^k(u)\bigr)
    \\
    &= \log \bigl(1+\lambda_{2n}\ch_{2n}(u)\bigr)
  \end{align*}
  In degree $4n$, the left\nobreakdash-hand side is
  $(k^{2n}-1)(-1)^{n-1} (B_n/2n)\ch_{2n}(u)$, so $\lambda_{2n} =
  (-1)^{n-1} B_n/2n$ as claimed.
\end{proof}

The two maps displayed in \eqref{eq:25} takes value in the
$(k-1)(1-g)$\nobreakdash-component of $\setZ\subprime \times
BU\subprime$ for surface bundles of fiber genus $g$.
Theorem~\ref{thm:A} is a consequence of diagram~\eqref{eq:23},
composed with $\alpha_E$, evaluated on $s_{2n-1}\in
H^{4n-2}(BU\subprime;\setZ)$. 

We begin by comparing the cohomology classes obtained by evaluating the
two maps
\begin{equation}
  \label{eq:30}
  \begin{split}
    f\colon & B \xrightarrow{\:\omega\circ \alpha_E\:}
    Q(\spaceCP_+^\infty) \xrightarrow{\:Q(\bar{L})\:}
    Q(\setZ\times BU) \xrightarrow{\:E\:}
    \setZ\subprime \times BU\subprime
    \\
    g\colon & B \xrightarrow{\:\omega\circ \alpha_E\:}
    Q(\spaceCP_+^\infty) \xrightarrow{\:Q(\bar{L})\:}
    Q(\setZ\times BU) \xrightarrow{\:Q(r^k)\:}
    Q(\setZ\subprime \times BU\subprime) \xrightarrow{\:E\:} 
    \setZ\subprime \times BU\subprime
  \end{split}
\end{equation}
on $s_{2n-1}\in H^*(BU\subprime;\setZ)$. By \eqref{eq:18} the first
map represents $\pi_* \bigl(\overline{\Tpi E} -1\bigr)$, the second is
one of the two maps in the homotopy commutating diagram \eqref{eq:23}.

The middle map $Q(\bar{L})$ is the sum of two:
\begin{align*}
  Q(\bar{L}-1)\colon & Q(\spaceCP_+^\infty) \to Q(BU),
  \\
  Q(1)\colon & Q(\spaceCP_+^\infty) \to Q(\setZ\times BU),
\end{align*}
and $Q(1)$ factors over the projection of $Q(\spaceCP_+^\infty)$ on
$Q(S^0)$. The composite
\begin{equation*}
  u\colon Q(S^0) \to Q(\setZ\times BU) \xrightarrow{\:E\:} \setZ\times
  BU 
\end{equation*}
is the infinite loop map of the connective unit in the
$K$\nobreakdash-theory ring spectrum $ku$. Similarly $E\circ
Q\bigl(r^k(1)\bigr)$ is the projection onto $Q(S^0))$ composed with
$u_k=\frac{k-1}{2} u$.

\begin{lemma}\label{lem:6.5}
  For odd primes $p$, $u^*(s_{2n-1})=0$ in
  $H^*\bigl(Q(S^0);\setZ\subprime\bigr)$. For $p=2$,
  $2u^*(s_{2n-1})=0$ in $H^*\bigl(Q(S^0);\setZ\subprime[2]\bigr)$ and
  $\omega^*u^*(s_{2n-1})=0$ in $H^*\bigl(\Omega^\infty MT(2);\setZ\subprime[2]\bigr)$.
\end{lemma}
The proof is deferred to the end of the section. The lemma implies
that in \eqref{eq:30} we may replace $Q(\bar{L})$ with $Q(\bar{L}-1)
\colon Q(\spaceCP_+^\infty)\to Q(BU)$ without affecting the evaluation
on $s_{2n-1}$. With this replacement, Lemma~\ref{lemma:key} and
Lemma~\ref{lem:6.3} combine to give
\begin{equation*}
  g^*(s_{2n-1})= (-1)^{n-1}(k^{2n}-1)B_n/2n\cdot f^*(s_{2n-1}).
\end{equation*}
The other map in \eqref{eq:23} is
\begin{equation*}
  h\colon B \xrightarrow{\:\alpha_E\:}
  \Omega^\infty MT(2) \xrightarrow{\:\lambda_{-L}\:}
  \setZ\times BU \xrightarrow{\:k\psi^k-\id\:}
  \setZ \times BU,
\end{equation*}
where, by \eqref{eq:25}, $\lambda_{-L}\circ \alpha_E$ represents
$\pi_*(1)$. Since $(k\psi^k-\id)^*$ multiplies $s_{2n-1}$ by
$(k^{2n}-1)$,
\begin{equation*}
  h^*(s_{2n-1}) = (k^{2n}-1)s_{2n-1}\bigl(\pi_*(1)\bigr),
\end{equation*}
and by \eqref{eq:23} $h^*(s_{2n-1})=g^*(s_{2n-1})$. This gives the
equation
\begin{equation}
  \label{eq:31}
  (k^{2n}-1)s_{2n-1}\bigl(\pi_*(1)\bigr)=(-1)^{n-1}(k^{2n}-1)B_n/2n
  \cdot s_{2n-1}\bigl(\pi_*(\overline{\Tpi E}-1)\bigr)
\end{equation}
in $H^*(B;\setZ\subprime)$.

\begin{proofof}[Proof of Theorem~\ref{thm:A}.]
  For $p$ odd, pick $k\in \setZ$ to represent a generator of the
  multiplicative group of units $(\setZ/p^2)^*$. Then
  \begin{equation*}
    (k^{2n}-1) = \Denom(B_n/2n)\cdot \lambda, \quad \lambda\in
    \setZ\subprime^* 
  \end{equation*}
  and \eqref{eq:31} gives
  \begin{equation*}
    \Denom(B_n/2n)s_{2n-1} \bigl(\pi_*(1)\bigr) = (-1)^{n-1}
    \Num(B_n/2n)s_{2n-1}\bigl(\pi_*(\overline{\Tpi E} -1)\bigr). 
  \end{equation*}
  If $p=2$, pick $k=\pm 3\pmod{8}$. Then
  \begin{equation*}
    (k^{2n}-1)=2\Denom(B_n/2n)\lambda
  \end{equation*}
  for a $2$\nobreakdash-local unit $\lambda$ and we get an extra
  factor of~2.  Since
  \begin{equation*}
    -s_{2n-1}\bigl(\pi_*(1)\bigr) = s_{2n-1}\bigl(\Hodge(E)\bigr)=: s_{2n-1}(E)
  \end{equation*}
  by Lemma~\ref{lem:3.4}, and
  \begin{equation*}
    \bStdClass_n(E) = s_n\bigl(\pi_*(\overline{\Tpi E}-1)\bigr) =
    s_n(\Index \bpartial _{\overline{\Tpi[] E}-1}).
  \end{equation*}
  This completes the proof.
\end{proofof}

We close the section with
\begin{proofof}[Proof of Lemma~\ref{lem:6.5}.]
  Pick $k$ to generate $(\setZ/p^2)^*$ when $p$ is odd, and $k=\pm
  3\pmod{8}$ for $p=2$. Define spaces $\spaceJ(p)$ and
  $\spaceJ[C](p)$ to make
  \begin{align*}
    \spaceJ(p) &\xrightarrow{\:i_{\setR}\:} 
    BO\subprime \xrightarrow{\:\psi^k-\id\:} B\Spin\subprime
    \\
    \spaceJ[C](p) & \xrightarrow{\:i_{\setC}\:}
    BU\subprime \xrightarrow{\:\psi^k-\id\:}
    BU\subprime
  \end{align*}
  homotopy fibrations. The space $Q(S^0)\subprime$ decomposes as a product
  \begin{equation}\label{eq:41}
    Q(S^0)\subprime \simeq \bigl(\setZ\subprime \times
    \spaceJ(p)\bigr) \times \cokernel \spaceJ(p)
  \end{equation}
  and $\tilde{K}\bigl(\cokernel\spaceJ(p)\bigr)=0$, see e.g.\
  \cite{MR0494076}, \cite{MR0494077}. Complexification induces a map
  $c\colon \spaceJ(p)\to\spaceJ[C](p)$, and there is a homotopy
  commutative diagram
  \begin{equation*}
    \xymatrix{
      Q(S^0) \ar[r] \ar[d]^{u} 
      &
      \setZ\subprime \times \spaceJ(p) \ar[d]^{c}
      \\
      \setZ\subprime\times BU \subprime 
      &
      \setZ\subprime \times \spaceJ[C](p) \ar[l]_-{i_{\setC}}
      \rlap{\quad .}
    }
  \end{equation*}
  The endormorphism $(\psi^k-\id)^*$ of $H^*(BU\subprime;\setZ)$
  multiplies the primitive generator $s_{2n-1}$ by $(k^{2n-1}-\nobreak
  1)$.  This is a unit of $\setZ\subprime$ if $p$ is odd, and twice a
  unit if $p=2$. Thus $2i_{\setC}^*(s_{2n-1})=0$ in all cases, and the
  diagram implies $2u^*(s_{2n-1})=0$.

  For $p=2$, we must further show that $\omega^*u^*(s_{2n-1})=0$. To
  this end, recall the fibration sequence
  \begin{equation*}
    \Omega^\infty MT(2) \xrightarrow{\ \omega\ } 
    Q(\spaceCP^\infty_+) \xrightarrow{\ \partial\ }
    \Omega Q_0(S^0),
  \end{equation*}
  where $\partial$ is the $S^1$\nobreakdash-transfer, \cite{ravenel},
  \cite{stoltz}. On the component $Q(S^0)$ of
  $Q(\spaceCP^\infty_+)=Q(\spaceCP^\infty)\times Q(S^0)$, $\partial$
  is the map $\eta\colon Q(S^0)\to\Omega Q_0(S^0)$ induced from the
  stable Hopf map $\eta\colon S^{n+1} \to S^n$. Applying the splitting
  \eqref{eq:41}, and the diagram
  \begin{equation*}
    \xymatrix{
      J_{\setR}(2) \ar[r] \ar[d]^\eta 
      &
      BO\subprime[2] \ar[r]^c \ar[d]^\eta
      &
      BU\subprime[2]
      \\
      \Omega J_{\setR} (2) \ar[r] 
      & 
      SO\subprime[2] \rlap{\quad ,}
    }
  \end{equation*}
  it suffices to check that the cohomology classes $c^*(s_{2n-1})\in
  H^*(BO\subprime[2];\setZ)$ lie in the image of $\eta\colon
  H^*(SO\subprime[2];\setZ)\to H^*(BO\subprime[2];\setZ)$. We can use
  the homotopy fibration sequence
  \begin{equation*}
    \cdots \longrightarrow BU \xrightarrow{\ r\ } BO \xrightarrow{\
      \eta\ } SO \xrightarrow{\ c\ } SU \longrightarrow \cdots
  \end{equation*}
  (which induces the Bott-sequence
  \begin{equation*}
    \cdots \longrightarrow\widetilde{K}(X)\xrightarrow{\ r\ } \widetilde{KO}(X)
    \longrightarrow
    \widetilde{KO}^{-1}(X) \longrightarrow\cdots \qquad )
  \end{equation*}
  to examine this question.

  The Serre spectral sequence of the homotopy fibration
  $BU\xrightarrow{\ r\ } BO \xrightarrow{\ \eta\ } SO$ has vanishing
  differentials, so that
  \begin{equation}
    \label{eq:42}
    H^*(BO,\setF_2)  \cong H^*(BU,\setF_2)\otimes H^*(SO;\setF_2)
  \end{equation}
  Indeed $H^*(SO;\setF_2)=\setF_2[x_1,x_2,\dots]$, so that the
  $E^2$\nobreakdash-term has the same dimension in each degree as the
  abutment $H^*(BO,\setF_2)$; this leaves no room for differentials.

  In integral cohomology, both $H^*(BO;\setZ\subprime[2])$ and
  $H^*(SO;\setZ\subprime[2])$ contains only torsion of order~2, and the
  Bockstein exact sequence implies that the vertical maps below are injective,
  \begin{equation*}
    \xymatrix{
      H^*(SO;\setZ\subprime[2]) \ar[r]^{\eta^*} 
      \ar@{>->}[d] 
      & 
      H^*(BO;\setZ\subprime[2]) \ar@{>->}[d]
      \\
      H^*(SO;\setF_2) \ar[r]^{\eta^*} 
      & 
      H^*(BO;\setF_2) \rlap{\quad.}
    }
  \end{equation*}
  The generators $x_{2k-1}\in H^{2k-1}(SO;\setF_2)$ are primitive and
  from \eqref{eq:42} it follows that in odd dimensions
  \begin{equation*}
    \eta^*\colon \Prim H^*(SO;\setF_2) \xrightarrow{\ \cong\ }
    \Prim H^*(BO;\setF_2)
  \end{equation*}
  is an isomorphism. Hence with $\setF_2$ coefficients
  \begin{equation*}
    \eta^*(x^2_{2n-1}) =c^*(s_{2n-1}) \text{ in $H^*(BO;\setF_2)$}. 
  \end{equation*}
  Since $x^2_{2n-1}$ is the image of a class in
  $H^*(SO;\setZ\subprime[2])$, this completes the proof.
\end{proofof}

\section{Discussion of Akita's conjecture} 
\label{sec:disc-akit-conj}

For a surface bundle $\pi\colon E\to B$, we have the standard
MMM\nobreakdash-classes
\begin{equation*}
  \StdClass_n(E) = \pi_*\bigl(c_1(\Tpi E)^{n+1}\bigr)\in H^{2n}(B;\setZ)
\end{equation*}
where
\begin{equation*}
  \pi_*\colon H^*(E;\setZ) \to H^{*-2}(B;\setZ)
\end{equation*}
is the cohomological push forward homomorphism. This section discusses
the relationship between these classes and the classes
$\bStdClass_n(E)$ that appear in Theorem~\ref{thm:A}.

Given a spectrum $A=\{A_n,\varepsilon_n\}$ with associated infinite
loop space $\Omega^\infty\! A$, one has the cohomology suspension map
\begin{equation*}
  \sigma^*\colon H^{k}(A;\setZ) \to H^{k}(\Omega^\infty\! A ;\setZ) .
\end{equation*}
The source is the spectrum cohomology,
\begin{equation*}
  H^k(A;\setZ)=\varprojlim H^{k+n}(A_n).
\end{equation*}
and $\sigma^*$ is induced from the evaluation $S^n\Omega^n\! A_n \to
A_n$. Dually, the homology suspension is the homomorphism
\begin{equation*}
  \sigma_*\colon H_k(\Omega^\infty\!A;\setZ)\to H_k(A;\setZ).
\end{equation*}
In terms of the cohomology suspension, the standard classes have the
following description
\begin{equation}
  \label{eq:32}
  \StdClass_n(E)=(\omega\circ \alpha_E)^*\sigma^*(e^n)
\end{equation}
where
\begin{equation*}
  B\xrightarrow{\:\alpha_E\:} \Omega^\infty \xrightarrow{\:\omega\:},
  Q(\spaceCP_+^\infty) 
\end{equation*}
and $e=c_1(L)$ is the Euler class of the tautological line bundle on
$\spaceCP^\infty$, cf.\ \cite[Thm.~3.1]{MR2219303}.

Fix an odd prime $p$. We shall compare $\StdClass_n(E)$ and
$\bStdClass_n(E)$ first in $H^*(B;\setF_p)$ and then in
$H^*(B;\setZ\subprime)$. We remember from
sect.~\ref{sec:proof-theor-refthm:a} that
\begin{equation}
  \label{eq:33}
  \bStdClass_n(E)=(\omega\circ \alpha_E)^*f^*(s_n),
\end{equation}
where $f\colon Q(\spaceCP_+^\infty)\xrightarrow{\:\proj\:}
Q(\spaceCP^\infty) \xrightarrow{Q(\bar{L}-1)\:} Q(BU)
\xrightarrow{\:E\:} BU$. 

\begin{lemma}\label{lem:7.3}
  In $H^*(B;\setQ)$, $\StdClass_n(E)=(-1)^n\bStdClass_n(E)$.
\end{lemma}
\begin{proof}
  There is a homotopy commutative diagram
  \begin{equation*}
    \xymatrix@C=2cm{
      Q(\spaceCP_+^\infty) \ar[r]^-{Q(\bar{L}-1)}
      &
      Q(BU) \ar[r]^-{E}
      & 
      BU
      \\
      \spaceCP^\infty \ar[u]_{i} \ar[r]^{\bar{L}-1}
      &
      BU \ar@{=}[ur] \ar[u]_{i}
      &
    }
  \end{equation*}
  with $i$ being the inclusions. Thus
  \begin{equation*}
    i^*Q(\bar{L}-1)E^*(s_n)=s_n(\bar{L})=(-1)^nc_1(L)^n\ .
  \end{equation*}
  But $i^*\sigma^* =\id$, so
  \begin{equation*}
    i^*\sigma^*(e^n)=e^n=c_1(L)^n.
  \end{equation*}
  Since $i^*\colon H^*\bigl(Q(X);\setQ\bigr))\to H^*(X;\setQ)$ is an
  isomorphism, the claim follows.
\end{proof}

We shall next consider the classes $\StdClass_n$ and $\bStdClass_n$ in
mod $p$ cohomology, $p$ odd. This requires infinite loop space theory,
or more precisely the theory of homology operations on the category of
infinite loop spaces. 
We restrict ourselves to odd primes $p$, leaving the case of $p=2$ for
the reader to work out.
These are homomorphisms
\begin{equation*}
  Q^i\colon H_k(\Omega^\infty\! A;\setF_p)\to
  H_{k+2i(p-1)}(\Omega^\infty\! A;\setF_p),
\end{equation*}
natural with respect to infinite loop maps (i.e.\ maps induced from
maps of spectra).

For any (pointed) space $X$, the
mod $p$ homology of $Q(X)$ is a functor of the mod $p$
homology of $X$ via iterated homology operations and Bocksteins,
\begin{equation*}
  Q^I=\beta^{\varepsilon_1}Q^{i_1}\beta^{\varepsilon_2}Q^{i_2}\cdots
  \beta^{\varepsilon_k}Q^{i_k}
  \qquad (\varepsilon_j = \text{$0$ or $1$, $i_j\in \setN$}),
\end{equation*}
applied to the classes of $H_*(X;\setF_p)\subset
H_*\bigl(Q(X);\setF_p\bigr)$, cf.\ \cite[I]{MR0436146}.

Assuming $X$ is connected, $H_*\bigl(Q(X);\setF_p\bigr)$ is the free
associative, graded commutative algebra on generators $Q^I(x)$ subject
to the following conditions for
$I=(\varepsilon_1,i_1,\varepsilon_2,i_2,\dots,\varepsilon_k,i_k)$:
\begin{equation}
  \label{eq:34}
  i_{j-1}\leq pi_j-\varepsilon_j,\qquad
  2s_1-\sum_{j=2}^k\bigl(2s_j(p-1)-\varepsilon_j\bigr)>\deg x,
\end{equation}
where $x\in \tilde{H}_*(X;\setF_p)\subset
\tilde{H}_*\bigl(Q(X);\setF_p\bigr)$ runs through a homogeneous vector
space basis. There is a similar statement for
non\nobreakdash-connected X, cf.\ \cite[I, Thm.~4.2]{MR0436146}.

It is well\nobreakdash-known that
\begin{equation*}
  \sigma_*\colon H_k(\Omega^\infty\!A;\setF_p)\to H_k(A;\setF_p)
\end{equation*}
annihilates products, and indeed all classes of the form $Q^i(x)$. The second
statement follows by iterating the homomorphism
\begin{equation*}
  \sigma_*\colon H_k(\Omega^\infty\!A;\setF_p) \to
  H_{k+1}(\Omega^{\infty-1}\! A;\setF_p).
\end{equation*}
Indeed $\sigma_*$ commutes with homology operations, and $Q^i(x)=0$
for $2i<\deg(x)$. The consequence we shall use is the following.

\begin{lemma}\label{lem:7.6}
  The classes $\sigma^*(e^n)\in
  H^{2n}\bigl(Q(\spaceCP_+^\infty);\setF_p\bigr)$ vanishes on all
    elements $Q^I(a_n)$ where $a_n\in H_{2n}(\spaceCP^\infty;\setF_p)$
      is the generator.
\end{lemma}
For infinite loop spaces $X=\Omega^\infty\!A$, the structure map
$E\colon Q(X)\to X$ is an infinite loop map so that
\begin{equation*}
  E_*\colon H_*\bigl(Q(X);\setF_p\bigr)\to H_*(X;\setF_p)
\end{equation*}
commutes with homology operations; the element $Q^I\bigl(i_*(x)\bigr)$
in $H_*\bigl(Q(X);\setF_p\bigr)$ is mapped to $Q^I(x)$. This will be
applied to $E\colon Q(BU)\to BU$ below.

Let $a_n\in H_{2n}(\spaceCP^\infty;\setZ)$ be dual to $e^n$,
$e=c_1(L)$, and denote also by $a_n$ its image under
\begin{equation*}
  \bar{L}-1\colon \spaceCP^\infty \to BU.
\end{equation*}
Then $H_*(BU;\setZ)$ is a polynomial algebra on the classes $a_n$.
Reducing to $\setF_p$ coefficients, \cite[Theorem~6]{MR0331386} proves the
formula:
\begin{equation}
  \label{eq:35}
  Q^j(a_n) = (-1)^{j+n-1}\tbinom{j-1}{n}a_{n+j(p-1)} +
\text{decomposables}, 
\end{equation}
where $\binom{a}{b}=\frac{a!}{b!(a-b)!}$ is the binomial coefficient
reduced to $\setF_p$.

We are now ready to compare the classes $\StdClass_n$ and $\bStdClass_n$ in
$H^{2n}\bigl(\Omega^\infty MT(2);\setZ)$. I will just do one concrete
example, but using the results from \cite{MR2079997} it is clear
that $\StdClass_n\neq \bStdClass_n$ for many values of~$n$.

\begin{proposition}\label{prop:7.8}
  Let $p$ be an odd prime.
  \begin{enumerate}
  \item \label{item:3} The classes $p\StdClass_n$ and $(-1)^np\bStdClass_n$ agree in
    $H^*\bigl(\Omega^\infty MT(2);\setZ\subprime\bigr)$.
  \item \label{item:4} $\StdClass_{2p-1}\neq -\bStdClass_{2p-1}$ in $H^*\bigl(\Omega^\infty
    MT(2);\setF_p)$.
  \end{enumerate}
\end{proposition}
\begin{proof}
  The vanishing of $p(\StdClass_n-(-1)^n\bStdClass_n)$ follows from
  Theorem~\ref{thm:5.5}. To  prove \ref{item:4}, we use the homotopy fibration
  \begin{equation*}
    \Omega^\infty MT(2) \xrightarrow{\:\omega\:} Q(\spaceCP_+^\infty)
    \xrightarrow{\partial_*} \Omega Q(S^0)\ ,
  \end{equation*}
  localised at the prime $p$. The space $\Omega Q(S^0)\subprime$ is
  $(2p-5)$\nobreakdash-connected, and $a_1\in
  H_2(\spaceCP^\infty;\setF_p) \subset
  H_2\bigl(Q(\spaceCP_+^\infty);\setF_p \bigr)$ is in the image of
  $\omega_*$, $a_1=\omega_*(\bar{a}_1)$. Then
  \begin{equation*}
    \innerp{\bStdClass_{2p-1},Q^2(\bar{a}_1)} = \innerp{s_{2p-1},Q^2(a_1)}
  \end{equation*}
  where the right\nobreakdash-hand side takes place in (co)homology of
  $BU$. By \eqref{eq:35},
  \begin{equation*}
    \innerp{\bStdClass_{2p-1},Q^2(\bar{a}_1)}
    =\innerp{s_{2p-1},a_{2p-1}} \neq 0\ .
  \end{equation*}
  On the other hand, $\innerp{\StdClass_{2p-1},Q^2(a_1)}=0$ since
  cohomology classes in the image of the suspension homomorphism
  annihilates homology operations.
\end{proof}

Note that the proposition above give a counterexample to Akita's conjecture,
stated in the introduction. This follows from \cite{MR2335797} where
it is proved that the map
\begin{equation*}
  \setZ\times B\groupDiff(F_g) \xrightarrow{\:\alpha\:} \Omega^\infty
  MT(2) 
\end{equation*}
is a homology isomorphism in degrees less than $(g-1)/2$. This implies
that there are surface bundles (of fiber genus at least $8p-3$) where
$\bStdClass_{2p-1} (E)\neq -\StdClass_{2p-1}(E)$.

\section[{Appendix on $p$\protect\nobreakdash-power operations}]{Appendix on $\bm{p}$\protect\nobreakdash-power operations}
\label{sec:append-pnobr-power}

This section is devoted to a proof of Theorem~\ref{thm:5.5}. The
properties enjoyed by $\PontryaginPpowerOp$ can be found in
\cite[sect.~1]{MR0375307} and \cite[sect.~I.7]{MR0436146}, which the
reader may consult for more details.

In this section $X$ will be a connected infinite loop space. The
$\setZ\subprime$\nobreakdash-module $H^*(X;\setZ\subprime)$ is in one
to one correspondence with the module of commutative diagrams
\begin{equation*}
  \xymatrix{
    H_*(X;\setQ) \ar[r] \ar[d]
    &
    \setQ \ar[d]
    \\
    H_*(X;\setZ/p^\infty) \ar[r]
    &
    \setZ/p^\infty
    \mathrlap{\ , \quad \setZ/p^\infty =\setQ/\setZ\subprime .   }
  }
\end{equation*}
The submodule of primitive elements $\Prim H^*(X;\setZ\subprime)$
corresponds to the diagrams
\begin{equation}
  \label{eq:36}
  \cxymatrix{
    QH_*(X;\setQ) \ar[r] \ar[d]
    &
    \setQ \ar[d]
    \\
    QH_*(X;\setZ/p^\infty) \ar[r]
    &
    \setZ/p^\infty
  }
\end{equation}
where $QH_*(X;?)$ denotes the module of indecomposable elements, cf.\
\cite[Lemma~1.3]{MR0375307}.

Let $E^r(X)$ denote the mod $p$ Bockstein spectral sequence associated
to the exact couple
\begin{equation*}
  \xymatrix@C=2mm{
    H_*(X;\setZ\subprime) \ar[rr]^-{p} 
    & &
    H_*(X;\setZ\subprime) \ar[ld]^{\varrho_1}
    \\
    & 
    H_*(X;\setZ/p)\ar[ul]_{\beta_1}
  }
\end{equation*}
where $\beta_1$ is the primary Bockstein operator. More generally, if
\begin{equation*}
  \beta\colon H_*(X;\setZ/p^\infty) \to H_*(X;\setZ\subprime)
\end{equation*}
is the universal Bockstein then $\beta_n = \beta\circ i_n$ with $i_n$
induced from $\setZ/p^n\subset \setZ/p^\infty$. Reduction mod $p$
defines a surjection of algebras
\begin{equation*}
  j^r\colon H_*(X;\setZ/p^r) \to E^r(X).
\end{equation*}
Its kernel is given by
\begin{equation}
  \label{eq:37}
  \Kernel(j^r) = \Image(p_*)+\Image(\varrho_r\circ \beta_{r-1}),
\end{equation}
where $\varrho_r$ is reduction to $\setZ/p^r$ coefficients and
$p_*\colon H_*(X;\setZ/p^{r-1})\to H_*(X;\setZ/p^r)$ is induced from
the inclusion of coefficients. There are commutative diagrams
\begin{equation}
  \label{eq:38}
  \cxymatrix{
    H_*(X;\setZ/p^r) \ar[r]^-{j^r} \ar[d]^{\varrho_r\circ \beta_r}
    & 
    E^r(X) \ar[d]^{d^r}
    \\
    H_*(X;\setZ/p^r) \ar[r]^-{j^r}
    &
    E^r(X)
  }
  \qquad\qquad
  \cxymatrix{
    H_*(X;\setZ/p^{r-1}) \ar[r]^-{j^{r-1}} \ar[d]^{p}
    & 
    E^{r-1}(X) \ar[d]^{\xi}
    \\
    H_*(X;\setZ/p^r) \ar[r]^-{j^r}
    &
    E^r(X)
  }
\end{equation}
($\xi(x)=x^p$).

Suppose now that $X=Q(Y)$ as in Theorem~\ref{thm:5.5}. The Bockstein
spectal sequence in this case is easy to describe, given the
description of $H_*(X;\setZ/p)$ in terms of homology operations
applied to the subspace $H_*(Y;\setZ/p)\subset H_*(X;\setZ/p)$:
$E^1(X)$ is a tensor product of two types of differential algebras, 
namely
\begin{enumerate}
\item \label{item:5} $P\{y_I\}\otimes E\{x_I\}$; $d^1y_I=x_I$ and
  $\deg y_I$ even,
\item \label{item:6} $E\{z_I\}\otimes P\{y_I\}$; $d^1z_I=y_I$ and
  $\deg y_I$ even.
\end{enumerate}
The first spectral sequence has $E^2=P\{y_I^p\}\otimes
E\{y_I^{p-1}x_I\}$, the second has $E^2=\setZ/p$, and
$d^2(y_I^p)=y_I^{p-1}x_I$. More generally, in case \ref{item:5}
\begin{equation}
  \label{eq:39}
  E^{r+1} = P\{y_I^{p^r}\} \otimes E\{y_I^{p^{r-1}}x_I\}, \quad 
  d^{r+1}(y_I^{p^r})=y_I^{p^r-1}x_I\ .
\end{equation}
There is a slight modification needed for $p=2$, cf.\
\cite[p.~49]{MR0436146} or \cite{MR0375307}, which we leave for the
reader. The Bockstein spectral sequence for $X=Q(Y)$ is from the
$E^2$\nobreakdash-term on a tensor product of the spectral sequences
\eqref{eq:39} with $y_I=Q^I(y)$,
\begin{equation}
  \label{eq:40}
  \begin{gathered}
    I=(i_1,\varepsilon_1,i_2,\varepsilon_2,\dots,\varepsilon_k,i_k),\quad
    i_{j-1}\leq pi_j-\varepsilon_j,
    \\
    2i_1-\sum_{j=2}^k\bigl(2i_j(p-1)-\varepsilon_j\bigr)>\deg y,
    \\
    \deg Q^I(y)\equiv 0 \pmod{2},
  \end{gathered}
\end{equation}
and with $y$ running over a vector space basis for
$\tilde{H}_*(Y;\setZ/p)$. This description of $E^r(X)$, $X=Q(Y)$ follows
easily from the formula for $H_*(X;\setZ/p)$ given in
sect.~\ref{sec:disc-akit-conj}, see also \cite[I,
Thm.~4.3]{MR0436146}. 

\begin{proofof}[Proof of Theorem~\ref{thm:5.5}]
  Let $X=Q(Y)$ with $Y$ connected and of finite type, and let $\gamma$ be a
  primitive element of $H_*(X;\setZ\subprime)$ which satisfies the
  three conditions \ref{item:7}, \ref{item:8} and \ref{item:9} of
  Theorem~\ref{thm:5.5}. The class $\gamma$ corresponds to a diagram
  \begin{equation*}
    \xymatrix{
      QH_*(X;\setQ) \ar[r]^-{\gamma_{\#}} \ar[d]
      & 
      \setQ \ar[d]
      \\
      QH_*(X;\setZ/p^\infty) \ar[r]^-{\gamma_{\#}}
      &
      \setZ/p^\infty \mathrlap{\ .}
    }
  \end{equation*}
  We remember that the inclusion $i\colon Y\to Q(Y)=X$ induces
  isomorphism in rational homology. As condition \ref{item:8} asserts
  that $i^*(\gamma)=0$,
  \begin{equation*}
    \gamma_{\#}\colon QH_*(X;\setQ)\to \setQ
  \end{equation*}
  is the zero homomorphism. It remains to show that the lower
  horizontal $\gamma_{\#}$ vanishes. Assume inductively that
  \begin{equation*}
    \gamma_{\#}\superprime[r]\colon H_*(X;\setZ/p^r)\to \setZ/p^r
  \end{equation*}
  is zero. This is true for $r=1$ by condition \ref{item:7}. We must
  verify it for $r$ replaced by $r+1$. Consider
  \begin{equation*}
    j^{r+1}\colon H_*(X;\setZ/p^{r+1}) \to E^{r+1}(X).
  \end{equation*}
  Its kernel is given in \eqref{eq:37} and the inductive assumption
  shows that $\gamma_{\smash{\#}}\superprime[r+1]$ also vanishes on
  $\Image(\varrho_{r+1}\circ\nobreak \beta_r)$, so
  $\gamma_{\#}\superprime[r+1]$ factor over $E^{r+1}(X)$, in fact over
  $QE^{r+1}(X)$ since $\gamma$ was assumed primitive. By
  \eqref{eq:39}, and \eqref{eq:40}
  \begin{equation*}
    QE^{r+1}(X) = \bigoplus_{I,y} \setZ/p \innerp{y_I^{p^r}} \oplus
    \setZ/p\innerp{y_I^{p^r-1}x_I}\ .
  \end{equation*}
  The diagrams in \eqref{eq:38} show that
  \begin{equation*}
    j^{r+1}\bigl(\PontryaginPpowerOp\superprime(y_I)\bigr) 
    = y_I^{p^r}, \quad
    j^{r+1}\bigl(\varrho_{r+1}\beta_{r+1}
    \PontryaginPpowerOp\superprime
    (y_I)\bigr)=y_I^{p^r-1}x_I \ .
  \end{equation*}
  By condition \ref{item:9},
  $\gamma_{\#}\bigl(\PontryaginPpowerOp\superprime(y_I)\bigr)=0$ and 
  $\gamma_{\#}$ annihilates
  $\varrho_{r+1}\beta_{r+1}\PontryaginPpowerOp\superprime(y_I)$ since
  $\gamma$ was an 
  integral class.
\end{proofof}

\newpage

\nocite{*}
\bibliographystyle{abbrv}
\bibliography{imadsen}

\begin{thebibliography}{10}

\bibitem{MR0198468}
J.~F. Adams.
\newblock On the groups {$J(X)$}. {II}.
\newblock {\em Topology}, 3:137--171, 1965.

\bibitem{MR505692}
J.~F. Adams.
\newblock {\em Infinite loop spaces}, volume~90 of {\em Annals of Mathematics
  Studies}.
\newblock Princeton University Press, Princeton, N.J., 1978.

\bibitem{MR1324104}
J.~F. Adams.
\newblock {\em Stable homotopy and generalised homology}.
\newblock Chicago Lectures in Mathematics. University of Chicago Press,
  Chicago, IL, 1995.
\newblock Reprint of the 1974 original.

\bibitem{MR1892095}
T.~Akita.
\newblock Nilpotency and triviality of mod {$p$} {M}orita-{M}umford classes of
  mapping class groups of surfaces.
\newblock {\em Nagoya Math. J.}, 165:1--22, 2002.

\bibitem{MR2405898}
T.~Akita and N.~Kawazumi.
\newblock Integral {R}iemann-{R}och formulae for cyclic subgroups of mapping
  class groups.
\newblock {\em Math. Proc. Cambridge Philos. Soc.}, 144(2):411--421, 2008.

\bibitem{MR1043170}
M.~F. Atiyah.
\newblock {\em {$K$}-theory}.
\newblock Advanced Book Classics. Addison-Wesley Publishing Company Advanced
  Book Program, Redwood City, CA, second edition, 1989.
\newblock Notes by D. W. Anderson.

\bibitem{MR0236950}
M.~F. Atiyah and I.~M. Singer.
\newblock The index of elliptic operators. {I}.
\newblock {\em Ann. of Math. (2)}, 87:484--530, 1968.

\bibitem{MR0236952}
M.~F. Atiyah and I.~M. Singer.
\newblock The index of elliptic operators. {III}.
\newblock {\em Ann. of Math. (2)}, 87:546--604, 1968.

\bibitem{MR0279833}
M.~F. Atiyah and I.~M. Singer.
\newblock The index of elliptic operators. {IV}.
\newblock {\em Ann. of Math. (2)}, 93:119--138, 1971.

\bibitem{MR0377873}
J.~C. Becker and D.~H. Gottlieb.
\newblock The transfer map and fiber bundles.
\newblock {\em Topology}, 14:1--12, 1975.

\bibitem{MR0420609}
J.~M. Boardman and R.~M. Vogt.
\newblock {\em Homotopy invariant algebraic structures on topological spaces}.
\newblock Lecture Notes in Mathematics, Vol. 347. Springer-Verlag, Berlin,
  1973.

\bibitem{MR0387496}
A.~Borel.
\newblock Stable real cohomology of arithmetic groups.
\newblock {\em Ann. Sci. \'Ecole Norm. Sup. (4)}, 7:235--272 (1975), 1974.

\bibitem{MR0436146}
F.~R. Cohen, T.~J. Lada, and J.~P. May.
\newblock {\em The homology of iterated loop spaces}.
\newblock Lecture Notes in Mathematics, Vol. 533. Springer-Verlag, Berlin,
  1976.

\bibitem{jebert06}
J.~Ebert.
\newblock {\em Characteristic classes of spin surface bundles: Applications of
  the Madsen-Weiss theory}.
\newblock PhD thesis, Universit{\"a}t Bonn, 2006.

\bibitem{MR648106}
O.~Forster.
\newblock {\em Lectures on {R}iemann surfaces}, volume~81 of {\em Graduate
  Texts in Mathematics}.
\newblock Springer-Verlag, New York, 1981.
\newblock Translated from the German by Bruce Gilligan.

\bibitem{MR2079997}
S.~Galatius.
\newblock Mod {$p$} homology of the stable mapping class group.
\newblock {\em Topology}, 43(5):1105--1132, 2004.

\bibitem{MR2219303}
S.~Galatius, I.~Madsen, and U.~Tillmann.
\newblock Divisibility of the stable {M}iller-{M}orita-{M}umford classes.
\newblock {\em J. Amer. Math. Soc.}, 19(4):759--779 (electronic), 2006.

\bibitem{art18}
S.~Galatius, I.~Madsen, U.~Tillmann, and M.~Weiss.
\newblock The homotopy type of the cobordism category, 2008.
\newblock To appear in Acta Mathematica.

\bibitem{MR1070716}
K.~Ireland and M.~Rosen.
\newblock {\em A classical introduction to modern number theory}, volume~84 of
  {\em Graduate Texts in Mathematics}.
\newblock Springer-Verlag, New York, second edition, 1990.

\bibitem{MR0331386}
S.~O. Kochman.
\newblock Homology of the classical groups over the {D}yer-{L}ashof algebra.
\newblock {\em Trans. Amer. Math. Soc.}, 185:83--136, 1973.

\bibitem{MR0375316}
H.~Ligaard and I.~Madsen.
\newblock Homology operations in the {E}ilenberg-{M}oore spectral sequence.
\newblock {\em Math. Z.}, 143:45--54, 1975.

\bibitem{MR0375307}
I.~Madsen.
\newblock Higher torsion in {$SG$} and {$BSG$}.
\newblock {\em Math. Z.}, 143:55--80, 1975.

\bibitem{MR1778114}
I.~Madsen and C.~Schlichtkrull.
\newblock The circle transfer and {$K$}-theory.
\newblock In {\em Geometry and topology: Aarhus (1998)}, volume 258 of {\em
  Contemp. Math.}, pages 307--328. Amer. Math. Soc., Providence, RI, 2000.

\bibitem{MR0494076}
I.~Madsen, V.~Snaith, and J.~Tornehave.
\newblock Infinite loop maps in geometric topology.
\newblock {\em Math. Proc. Cambridge Philos. Soc.}, 81(3):399--430, 1977.

\bibitem{MR2335797}
I.~Madsen and M.~Weiss.
\newblock The stable moduli space of {R}iemann surfaces: {M}umford's
  conjecture.
\newblock {\em Ann. of Math. (2)}, 165(3):843--941, 2007.

\bibitem{MR0420610}
J.~P. May.
\newblock {\em The geometry of iterated loop spaces}.
\newblock Springer-Verlag, Berlin, 1972.
\newblock Lectures Notes in Mathematics, Vol. 271.

\bibitem{MR0494077}
J.~P. May.
\newblock {\em {$E\sb{\infty }$} ring spaces and {$E\sb{\infty }$} ring
  spectra}.
\newblock Lecture Notes in Mathematics, Vol. 577. Springer-Verlag, Berlin,
  1977.
\newblock With contributions by Frank Quinn, Nigel Ray, and J\o rgen Tornehave.

\bibitem{MR914849}
S.~Morita.
\newblock Characteristic classes of surface bundles.
\newblock {\em Invent. Math.}, 90(3):551--577, 1987.

\bibitem{MR717614}
D.~Mumford.
\newblock Towards an enumerative geometry of the moduli space of curves.
\newblock In {\em Arithmetic and geometry, Vol. II}, volume~36 of {\em Progr.
  Math.}, pages 271--328. Birkh\"auser Boston, Boston, MA, 1983.

\bibitem{ravenel}
D.~C. Ravenel.
\newblock The {S}egal conjecture for cyclic groups and its consequences.
\newblock {\em Amer. J. Math.}, 106(2):415--446, 1984.
\newblock With an appendix by Haynes R. Miller.

\bibitem{MR0353298}
G.~Segal.
\newblock Categories and cohomology theories.
\newblock {\em Topology}, 13:293--312, 1974.

\bibitem{stoltz}
S.~Stolz.
\newblock {\em Beziehungen zwischen {T}ransfer und {$J$}-{H}omomorphismus}.
\newblock Bonner Mathematische Schriften [Bonn Mathematical Publications], 125.
  Universit\"at Bonn Mathematisches Institut, Bonn, 1980.
\newblock Dissertation, Universit{\"a}t Bonn, Bonn, 1979.

\bibitem{MR1000393}
T.~tom Dieck.
\newblock The {A}rtin-{H}asse logarithm for {$\lambda$}-rings.
\newblock In {\em Algebraic topology (Arcata, CA, 1986)}, volume 1370 of {\em
  Lecture Notes in Math.}, pages 409--415. Springer, Berlin, 1989.

\end{thebibliography}

\noindent
\textsc{University of Copenhagen, Copenhagen, Denmark}

\noindent
\textsc{E-mail address:} imadsen@math.ku.dk

\end{document}